# ESTIMATION OF SUMS OF RANDOM VARIABLES: EXAMPLES AND INFORMATION BOUNDS[1]

By Cun-Hui Zhang

*Rutgers University*

This paper concerns the estimation of sums of functions of observable and unobservable variables. Lower bounds for the asymptotic variance and a convolution theorem are derived in general finite- and infinite-dimensional models. An explicit relationship is established between efficient influence functions for the estimation of sums of variables and the estimation of their means. Certain "plug-in" estimators are proved to be asymptotically efficient in finite-dimensional models, while "$u,v$" estimators of Robbins are proved to be efficient in infinite-dimensional mixture models. Examples include certain species, network and data confidentiality problems.

**1. Introduction.** Given a pool of $n$ motorists, how do we estimate the total intensity of those in the pool who have a prespecified number of traffic accidents in a given time period? This is an example of a broad class of problems involving the estimation of sums of random variables

$$(1.1) \qquad S_n \equiv \sum_{j=1}^n u(X_j, \theta_j)$$

[24], where $X_j$ are observable variables, $\theta_j$ are unobservable variables or constants, and $u(\cdot, \cdot)$ is a certain utility function. The estimation of (1.1) has numerous important applications. In the motorist example, $X_j$ is the number of traffic accidents and $\theta_j$ the intensity of the $j$th individual in the pool, and $u(x, \vartheta) = \vartheta I\{x = a\}$ for a prespecified integer $a$. In Sections 3, 4 and 5 we consider applications in certain species, network and data confidentiality problems.

Received June 2001; revised October 2004.
[1]Supported in part by the National Science Foundation.

The estimation of (1.1) is a nonstandard problem in statistics, since the sums, involving observables, as well as unobservables, are not parameters. Without a theory of efficient estimation, the performance of different estimators can only be measured against each other in terms of relative efficiency. For the specific motorist example with $u(x,\vartheta) = \vartheta I\{x = a\}$, Robbins and Zhang [28] proved that, in a Poisson mixture model, the efficient estimation of (1.1) is equivalent to the efficient estimation of $E(\theta|X = a)$, so that the usual information bounds can be used. In this paper we provide a general theory for the efficient estimation of sums of variables.

Let $(X, \theta)$, $(X_j, \theta_j)$, $j = 1, \ldots, n$, be i.i.d. vectors with an unknown common joint distribution $F$. Our general theory covers asymptotic efficiency for the estimation of

$$S_n \equiv S_n(F) \equiv \sum_{j=1}^{n} u(X_j, \theta_j; F) \tag{1.2}$$

based on $X_1, \ldots, X_n$, where the utility $u(x, \vartheta; F)$ is also allowed to depend on $F$. This provides a unified asymptotic theory for the estimation of (1.1) and conventional parameters $u(F)$, since the utility is allowed to depend on $F$ only. Our problem is closely related to the estimation of the mean

$$\mu(F) \equiv E_F u(X, \theta; F). \tag{1.3}$$

If $E_F u^2(X, \theta; F) < \infty$ and $1/2 \leq \alpha < 1$, an estimator is $n^\alpha$-consistent for the estimation of $S_n(F)$ iff it is $n^\alpha$-consistent for the estimation of its mean $n\mu(F) = E_F S_n(F)$. But an efficient estimator of $n\mu(F)$ is not necessarily an efficient estimator of $S_n(F)$, since the two estimation problems may have different efficient influence functions, as we demonstrate below in (1.4)–(1.6) and in simple examples in Sections 2.3 and 2.4. The asymptotic theory for the estimation of $\mu(F)$ is well understood; see [3, 17, 31].

Suppose that $F$ belongs to a known class $\mathcal{F}$. Let $F_0 \in \mathcal{F}$. An estimator $\widehat{\mu}_n$ of (1.3) is (locally) asymptotically efficient in contiguous neighborhoods of $P_{F_0}$ iff

$$\widehat{\mu}_n = \mu(F_0) + \frac{1}{n} \sum_{j=1}^{n} \psi_*(X_j) + o_{P_{F_0}}(n^{-1/2}), \tag{1.4}$$

where $\psi_*(x) \equiv \psi_*(x; F_0)$ is the efficient influence function at $F_0$ for the estimation of $\mu(F)$. In Section 6 we show that, under mild regularity conditions on the utility functions $\{u(x, \vartheta; F), F \in \mathcal{F}\}$, an estimator $\widehat{S}_n$ of (1.2) is (locally) asymptotically efficient in contiguous neighborhoods of $P_{F_0}$ iff

$$\frac{\widehat{S}_n}{n} = \mu(F_0) + \frac{1}{n} \sum_{j=1}^{n} \phi_*(X_j) + o_{P_{F_0}}(n^{-1/2}), \tag{1.5}$$

where $\phi_*(x) \equiv \phi_*(x; F_0)$ is the efficient influence function at $F_0$ for the estimation of $S_n(F)$. Furthermore, the following relationship holds between the two efficient influence functions in (1.4) and (1.5):

$$\phi_*(x) = \psi_*(x) + \overline{u}(x; F_0) - \mu(F_0) - u_*(x), \tag{1.6}$$

where $\overline{u}(x; F) \equiv E_F[u(X, \theta; F)|X = x]$ and $u_*(x) \equiv u_*(x; F_0)$ is the projection of $\overline{u}(x; F_0)$ to the tangent space of the family of distributions $\{F^X, F \in \mathcal{F}\}$ at $F_0^X$. Here $F^X$ is the marginal distribution of $X$ under the joint distribution $F$ of $(X, \theta)$. It follows clearly from (1.6) that asymptotically efficient estimations of $S_n(F)/n$ and $\mu(F)$ are equivalent in contiguous neighborhoods of $P_{F_0}$ iff $\overline{u}(\cdot; F_0) - \mu(F_0)$ is in the tangent space, that is, $\overline{u}(\cdot; F_0) - \mu(F_0) = u_*(\cdot; F_0)$.

We will derive more explicit results in finite-dimensional models and infinite-dimensional mixture models. In finite-dimensional models $\mathcal{F} = \{F_\tau, \tau \in \mathcal{T}\}$ with a Euclidean $\tau$, it will be shown that "plug-in" estimators of the form $\sum_{j=1}^n \overline{u}(X_j; F_{\widehat{\tau}_n})$ are asymptotically efficient for the estimation of (1.2) if $\widehat{\tau}_n$ is an efficient estimator of $\tau$. In infinite-dimensional mixture models, certain "$u, v$" estimators of Robbins [24] will be shown to be efficient for the estimation of (1.1). We shall consider estimation of (1.1) with known $f(x|\vartheta)$ in Section 2 and provide the general theory in Section 6. Section 7 contains proofs of all theorems.

**2. Mixture models.** Suppose $(X, \theta) \sim F(dx, d\vartheta) = f(x|\vartheta)\nu(dx)G(d\vartheta)$, that is,

$$X|\theta \sim f(x|\theta), \qquad \theta \sim G. \tag{2.1}$$

In this section we state our results for the estimation of (1.1) with known $f(\cdot|\cdot)$.

2.1. *Finite-dimensional mixture models.* Let $\{G_\tau, \tau \in \mathcal{T}\}$ be a parametric family of distributions with an open $\mathcal{T}$ in a Euclidean space. Suppose (2.1) holds with $G = G_\tau$ for an unknown vector $\tau \in \mathcal{T}$. Suppose that, for certain functions $\widetilde{\rho}_\tau$,

$$\int (\sqrt{g_{\tau, \Delta}} - 1 - \Delta^t \widetilde{\rho}_\tau/2)^2 \, dG_\tau = o(\|\Delta\|^2), \tag{2.2}$$
$$\int g_{\tau, \Delta} \, dG_\tau = 1 + o(\|\Delta\|^2), \qquad \text{as } \Delta \to 0,$$

where $g_{\tau, \Delta}$ is the Radon–Nikodym derivative of the absolutely continuous part of $G_{\tau+\Delta}$ with respect to $G_\tau$. Let $E_\tau$ denote the expectation under $G_\tau$. The Fisher information matrix for the estimation of $\tau$ based on a single $X$ is

$$I_\tau \equiv \operatorname{Cov}_\tau(\rho_\tau(X)), \qquad \rho_\tau(x) \equiv E_\tau[\widetilde{\rho}_\tau(\theta)|X = x]. \tag{2.3}$$

Define $\overline{u}_\tau(x) \equiv E_\tau[u(X, \theta)|X = x]$ and $\mu_\tau \equiv E_\tau u(X, \theta)$.

THEOREM 2.1. *Suppose* (2.2) *holds,* $E_\tau u^2(X,\theta)$ *is locally bounded and* $I_\tau$ *are of full rank for all* $\tau \in \mathcal{T}$. *Then* $\{\widehat{S}_n, n \geq 1\}$ *is an asymptotically efficient estimator of* (1.1) *iff* (1.5) *holds with* $\mu(F_0) = \mu_\tau$, $P = P_\tau$, *and the efficient influence function*

$$(2.4) \qquad \phi_* = \phi_{*,\tau} \equiv \overline{u}_\tau - \mu_\tau + \rho_\tau^t I_\tau^{-1} \gamma_\tau,$$

*where* $\gamma_\tau \equiv E_\tau \operatorname{Cov}_\tau(u(X,\theta), \widetilde{\rho}_\tau(\theta)|X) = E_\tau\{u(X,\theta)\widetilde{\rho}_\tau(\theta) - \overline{u}_\tau(X)\rho_\tau(X)\}$.

REMARK 2.1. Since $\kappa_{*,\tau} \equiv I_\tau^{-1}\rho_\tau$ is the efficient influence function for the estimation of $\tau$ and $\partial \mu_\tau/\partial \tau = E_\tau U(X,\theta)\widetilde{\rho}_\tau(\theta)$, $\psi_{*,\tau} \equiv \rho_\tau^t I_\tau^{-1} E_\tau u(X,\theta)\widetilde{\rho}_\tau(\theta)$ is the efficient influence function for the estimation of $\mu_\tau$. Moreover, $\overline{u}_{*,\tau} \equiv \rho_\tau^t I_\tau^{-1} E_\tau \overline{u}_\tau(X)\rho_\tau(X)$ is the projection of $\overline{u}_\tau$ to the tangent space generated by the scores $\rho_\tau(X)$ under $E_\tau$. Thus, Theorem 2.1 asserts that (1.5) and (1.6) hold under (2.2).

Our next theorem provides the asymptotic theory for plug-in estimators

$$(2.5) \qquad \widehat{S}_n \equiv \sum_{j=1}^{n} \overline{u}_{\widehat{\tau}_n}(X_j)$$

of (1.1), where $\overline{u}_\tau(x) \equiv E_\tau[u(X,\theta)|X=x]$ as in Theorem 2.1. An estimator $\widehat{\tau}_n$ of the vector $\tau$ is an asymptotically linear one with influence functions $\kappa_\tau$ under $E_\tau$ if

$$(2.6) \qquad \widehat{\tau}_n = \frac{1}{n}\sum_{j=1}^{n}\kappa_\tau(X_j) + o_{P_\tau}(n^{-1/2}),$$

with $E_\tau \kappa_\tau(X)\rho_\tau^t(X)$ being the identity matrix.

THEOREM 2.2. *Let* $\widehat{S}_n$ *be as in* (2.5) *with an asymptotically linear estimator* $\widehat{\tau}_n$ *as in* (2.6). *Suppose conditions of Theorem 2.1 hold,* $E_\tau \overline{u}_{\tau+\Delta}^2(X) = O(1)$ *as* $\Delta \to 0$ *for every* $\tau \in \mathcal{T}$, *and for all* $\tau \in \mathcal{T}$ *and* $c > 0$,

$$(2.7) \qquad \sup_{\|\Delta\|\leq c/\sqrt{n}}\left|\sum_{j=1}^{n}[\overline{u}_{\tau+\Delta}(X_j) - \overline{u}_\tau(X_j) - \{E_\tau\overline{u}_{\tau+\Delta}(X) - \mu_\tau\}]\right| = o_{P_\tau}(n^{1/2}).$$

*Let* $\phi_{*,\tau}$ *and* $\gamma_\tau$ *be as in Theorem 2.1 and* $\kappa_{*,\tau} = I_\tau^{-1}\rho_\tau$. *Then*

$$(2.8) \qquad \frac{\widehat{S}_n - S_n}{n^{1/2}} \xrightarrow{D} N(0,\sigma_\tau^2), \qquad \sigma_\tau^2 = \sigma_{*,\tau}^2 + \operatorname{Var}_\tau(\{\kappa_\tau(X) - \kappa_{*,\tau}(X)\}^t \gamma_\tau)$$

*under* $E_\tau$, *where* $\sigma_{*,\tau}^2 \equiv \operatorname{Var}_\tau(\phi_{*,\tau}(X) - u(X,\theta))$. *Consequently,* $\widehat{S}_n$ *is an asymptotically efficient estimator of* (1.1) *at* $E_{\tau_0}$ *iff* $\gamma_{\tau_0}\widehat{\tau}_n$ *is an asymptotically efficient estimator of* $\gamma_{\tau_0}\tau$ *in contiguous neighborhoods of* $E_{\tau_0}$.

REMARK 2.2. It follows from (2.8) that $|\widehat{S}_n - S_n| \leq 1.96 \sigma_{\widehat{\tau}_n} n^{1/2}$ provides an approximate 95% confidence interval for (1.1), provided that $\sigma_\tau$ is continuous in $\tau$.

REMARK 2.3. Condition (2.7) holds if $\{\overline{u}_{\tau+\Delta} : \tau + \Delta \in \mathcal{T}, \|\Delta\| \leq \delta_\tau\}$ is a Donsker class under $E_\tau$ for some $\delta_\tau > 0$ and $E_\tau \overline{u}_{\tau+\Delta}^2(X)$ is continuous at $\Delta = 0$.

2.2. *General mixtures.* Let $\mathcal{G}$ be a convex class of distributions. Suppose (2.1) holds with an unknown $G \in \mathcal{G}$. Let $E_G$ be the expectation under (2.1). Suppose $E_G u^2(X, \theta) < \infty$ for all $G \in \mathcal{G}$. Define

$$(2.9) \quad \mathcal{G}_{G_0} \equiv \left\{ G : E_{G_0}(f_G(X)/f_{G_0}(X))^2 < \infty, \int f_G I\{f_{G_0} > 0\} \, d\nu = 1 \right\},$$

where $f_G(x) \equiv \int f(x|\vartheta) G(d\vartheta)$, and define

$$(2.10) \quad \mathcal{V}_{G_0} \equiv \{v(x) : E_G v(X) = E_G u(X, \theta) \, \forall G \in \mathcal{G}_{G_0}\}.$$

THEOREM 2.3. (i) *If $\mathcal{V}_{G_0}$ is nonempty, then $\{\widehat{S}_n, n \geq 1\}$ is an asymptotically efficient estimator of (1.1) at $E_{G_0}$ iff $\widehat{S}_n = \{\sum_{j=1}^n v_{G_0}(X_j)\} + o_{P_{G_0}}(n^{1/2})$ with*

$$(2.11) \quad v_{G_0} \equiv \arg\min\{E_{G_0}(v(X) - u(X, \theta))^2 : v \in \mathcal{V}_{G_0}\}.$$

(ii) *If $\mathcal{V}_{G_0}$ is empty, then there does not exist any regular $n^{-1/2}$-consistent estimator of $E_G u(X, \theta)$ or $S_n/n$ in contiguous neighborhoods of $E_{G_0}$.*

The definition of regular estimators of (1.1) is given in Section 6.
Suppose that for certain $\mathcal{G}_* \subseteq \mathcal{G}$ the collection

$$(2.12) \quad \mathcal{V}_* \equiv \{v(x) : E_G v(X) = E_G u(X, \theta), E_G v^2(X) < \infty \, \forall G \in \mathcal{G}_*\}$$

is nonempty, for example, certain $\mathcal{V}_{G_0}$ as in Theorem 2.3(i). Let $\|h\|_G \equiv \{E_G h^2(X)\}^{1/2}$.

THEOREM 2.4. *Let $v_{G_0}$ be as in (2.11). Suppose $v_{G_0} \in \mathcal{V}_*$ and as $(\varepsilon, n) \to (0, \infty)$,*

$$\sup\left\{\left|\sum_{j=1}^n \frac{v_G(X_j) - v_{G_0}(X_j)}{n^{1/2}}\right| : \|v_G - v_{G_0}\|_{G_0} \leq \varepsilon, G \in \mathcal{G}_*\right\} \to 0 \quad \text{in } P_{G_0}$$

*for all $G_0 \in \mathcal{G}_*$. Let $\widehat{G}$ be an estimator of $G$ such that $P_{G_0}(\widehat{G} \in \mathcal{G}_*) \to 1$ and $\|v_{\widehat{G}} - v_{G_0}\|_{G_0} \to 0$ in $P_{G_0}$ for all $G_0 \in \mathcal{G}_*$. Then*

$$(2.13) \quad \widehat{V}_n \equiv \sum_{j=1}^n v_{\widehat{G}}(X_j)$$

*is an asymptotically efficient estimator of* (1.1) *at* $P_{G_0}$ *for all* $G_0 \in \mathcal{G}_*$.

If $f(x|\vartheta)$ belongs to certain exponential families, there exists a unique function $v$ such that $\mathcal{V}_{G_0} \neq \varnothing$ implies $\mathcal{V}_{G_0} = \{v\}$, so that $v_{G_0} = v$ for all $G_0$ and $\mathcal{V}_* = \{v\}$. The following theorem is a variation of Theorem 2.4 for such distributions.

THEOREM 2.5. *Suppose* $f(x|\vartheta) \propto \exp(x^t \lambda(\vartheta))$, $\lambda(\vartheta) \in \Lambda$, *is an exponential family with an open* $\Lambda$ *in a Euclidean space, and that the conditional distribution of* $\theta$ *given* $\lambda(\theta)$ *is known. Suppose* $\mathcal{G}$ *contains distributions* $G \equiv G_c$ *with* $E_G |\lambda(\theta) - c| = 0$ *for all* $c \in \Lambda$. *If* $\mathcal{V}_{G_0} \neq \varnothing$ *for certain* $G_0$, *then there exists a function* $v(x)$ *such that*

$$(2.14) \quad E_G[v(X)|\lambda(\theta) = c] = E_G[u(X,\theta)|\lambda(\theta) = c] \qquad \forall c \in \Lambda, G \in \mathcal{G},$$

*and such that the following* $V_n$ *is an efficient estimator of* $S_n$ *under* $\{E_G : E_G v^2(X) < \infty\}$:

$$(2.15) \qquad\qquad V_n \equiv \sum_{j=1}^n v(X_j).$$

REMARK 2.4. Robbins [24] called (2.15) "$u, v$" estimators, provided that (2.14) holds. The $\widehat{V}_n$ in (2.13) can be viewed as a "$u, v$" estimator with an estimated optimal $v$. Theorems 2.4 and 2.5 provide conditions under which these two types of "$u, v$" estimators are asymptotically efficient.

2.3. *The Poisson example.* Let $(X, Y, \lambda) \equiv (X, \theta)$ with

$$(2.16) \quad \begin{aligned} E[Y|X,\lambda] &= \lambda, \\ f(x|\lambda) &\equiv P(X = x|\lambda) = e^{-\lambda}\lambda^x/x!, \qquad x = 0, 1, \ldots. \end{aligned}$$

Robbins [22, 24] and Robbins and Zhang [25, 26, 27] considered the estimation of $S'_n \equiv \sum_{j=1}^n \lambda_j u(X_j)$ and $S''_n \equiv \sum_{j=1}^n Y_j u(X_j)$, and several related problems.

Both $S'_n$ and $S''_n$ are special cases of (1.1). For $u(x) = I\{x \leq a\}$, $S''_n$ could be the total number of accidents next year for those motorists with no more than $a$ accidents this year in the motorist example.

Suppose $\lambda_j$ have a common exponential density $\tau e^{-\lambda \tau} d\lambda$ with unknown $\tau$. The marginal distribution of $X$ is $f_\tau(x) = \tau(1+\tau)^{-x-1}$, and the marginal and conditional expectations of $\lambda u(X)$ and $Yu(X)$ are

$$\overline{u}_\tau(x) = \frac{(x+1)u(x)}{1+\tau}, \qquad \mu_\tau = \sum_{x=0}^\infty f_\tau(x) x u(x-1).$$

Let $\overline{X} \equiv \sum_{j=1}^{n} X_j/n$. Define $\widehat{\tau}_n \equiv (\beta+n)/(\alpha+\sum_{j=1}^{n} X_j)$ and

$$(2.17) \qquad \widehat{S}_n \equiv \sum_{j=1}^{n} \overline{u}_{\widehat{\tau}_n}(X_j) = \sum_{j=1}^{n} \frac{(\alpha/n+\overline{X})(X_j+1)u(X_j)}{(\alpha+\beta)/n+1+\overline{X}}.$$

It follows from Theorem 2.2 that the plug-in estimators in (2.17) are asymptotically efficient for both $S'_n$ and $S''_n$. For $\alpha = \beta = 0$, (2.17) gives the plug-in estimator corresponding to the maximum likelihood estimator (MLE) of $\tau$. For general positive $\alpha$ and $\beta$, (2.17) gives the Bayes estimator of $S'_n$ and $S''_n$ with a beta prior on $\tau/(1+\tau)$. Clearly, $\widehat{\mu}_n \equiv \sum_{x=1}^{\infty} \{\widehat{\tau}_n x u(x-1)\}/(1+\widehat{\tau}_n)^{x+1}$ is efficient for the estimation of the mean $\mu_\tau \equiv E_\tau u(X,\theta)$, but not for $S'_n/n$ or $S''_n/n$. Similar results can be obtained for $\lambda$ with the gamma distribution; see [23].

In the case of completely unknown $G(d\lambda)$, the "$u,v$" estimator (2.15) with $v(x) = xu(x-1)$ is asymptotically efficient for the estimation of $S'_n$ and $S''_n$ for all $G$ with finite $E_G\{v(X) - \lambda u(X)\}^2$.

### 2.4. More examples.

EXAMPLE 2.1. Let $X \sim N(\tau, \sigma^2)$. The number of "above average" individuals, $\widehat{S}_n \equiv \#\{j \leq n : X_j > \overline{X}\}$, is an efficient estimator of the number of above mean individuals $S_n(\tau) \equiv \#\{j \leq n : X_j > \tau\}$. The estimator $\widetilde{S}_n \equiv n/2$ is efficient for the estimation of $E_\tau S_n(\tau) = n/2$, but not $S_n(\tau)$.

EXAMPLE 2.2. Let $f(x|\vartheta) \sim N(\vartheta, \sigma^2)$. An efficient estimator for the number of "above mean" individuals, $S_n \equiv \#\{j \leq n : X_j > \theta_j\}$, is $\widehat{S}_n \equiv n/2$, compared with Example 2.1. This is even true under the condition $n^{-1}\sum_{j=1}^{n} \theta_j^2 = O(1)$, that is, in contiguous neighborhoods of $P_0$ with $P_0\{\theta_j = 0\} = 1$.

EXAMPLE 2.3. $\widehat{S}_n \equiv 0$ is efficient for the estimation of $S_n(\tau) \equiv \sum_{j=1}^{n} \rho_\tau(X_j)$.

**3. A species problem.** An interesting example of our problem is estimating the total number of species in a population of plants or animals. Suppose a random sample of size $N$ is drawn (with replacement) from a population of $d$ species. Let $n_k$ be the number of species represented $k$ times in the sample. A species problem is to estimate $d$ based on $\{n_k, k \geq 1\}$. The problem dates back to [13] and [14] and has many important applications [4]. We consider a network application in Section 4.

3.1. *Finite-dimensional models.* Let $X_j$ be the frequencies of the $j$th species in the sample, so that, for certain $p_j > 0$,

$$(3.1) \quad n_k = \sum_{j=1}^{d} I\{X_j = k\}, \qquad (X_1, \ldots, X_d) \sim \text{multinomial}(N, p_1, \ldots, p_d).$$

We will confine our discussion to the case of $(N, N/d) \to (\infty, \mu)$, $0 < \mu < \infty$, since $E(d - \sum_{k=1}^{\infty} n_k) = \sum_{j=1}^{d}(1-p_j)^N \to 0$ as $N \to \infty$ for fixed $d$. Let $\{G_\tau, \tau \in \mathcal{T}\}$ be a parametric family of distributions in $(0, \infty)$, where $\tau$ is an unknown parameter with a scale component, $G_\tau(y/c) = G_{\tau'_c}(y)$. Let $P_\tau$ be probability measures under which (3.1) holds conditionally on $N$ and certain i.i.d. variables $\theta_j > 0$, and

$$(3.2) \quad p_j = \frac{\theta_j}{\sum_{i=1}^{d} \theta_i}, \qquad N|\{\theta_j\} \sim \text{Poisson}\left(c \sum_{j=1}^{d} \theta_j\right), \qquad \theta_j \sim G,$$

with $G = G_\tau$. Under $P_\tau$, $X_j$ are i.i.d. with $P_\tau\{X_j = k\} = \int e^{-y}(y^k/k!)G_{\tau'_c}(dy)$. Assume $c = 1$ due to scale invariance. Since $n_0$ is unobservable, the MLE of $(d, \tau)$ is

$$(3.3) \quad \widehat{d} \equiv \frac{\sum_{k=1}^{N} n_k}{\int(1 - e^{-y})G_{\widehat{\tau}}(dy)}, \qquad \widehat{\tau} \equiv \arg\max_{\tau \in \mathcal{T}} \prod_{k=1}^{\infty} \left\{\frac{\int e^{-y} y^k G_\tau(dy)}{1 - \int e^{-y} G_\tau(dy)}\right\}^{n_k}.$$

In the next two paragraphs we derive the influence function for the MLE (3.3) and prove its asymptotic efficiency.

If (2.2) holds and the MLE $\widehat{\tau}$ of $\tau$ is asymptotically efficient, then

$$(3.4) \qquad \widehat{\tau} = \tau + \frac{1}{d}\sum_{i=1}^{d}\kappa_{*,\tau}(X_j) + o_P(d^{-1/2})$$

with $\kappa_{*,\tau} \equiv \{\text{Cov}_\tau(\overline{\rho}_\tau(X))\}^{-1}\overline{\rho}_\tau$ and $\overline{\rho}_\tau \equiv I_{\{x>0\}}(\rho_\tau(x) - \gamma_\tau)$, where $\rho_\tau$ is as in (2.3) and $\gamma_\tau \equiv E_\tau[\rho_\tau(X)|X > 0]$. Thus, by the Taylor expansion of the $\widehat{d}$ in (3.3),

$$(3.5) \qquad \widehat{d} = d + \sum_{j=1}^{d} \phi_{*,\tau}(X_j) + o_P(d^{1/2}),$$

where $\phi_{*,\tau}(x) \equiv I_{\{x>0\}}/P_\tau(X > 0) - 1 - \kappa_{*,\tau}^t(x)\gamma_\tau$. In this case, as $d \to \infty$,

$$(3.6) \qquad \frac{\widehat{d} - d}{d^{1/2}} \xrightarrow{D} N\left(0, \frac{P_\tau(X = 0)}{P_\tau(X > 0)} + \gamma_\tau^t\{\text{Cov}_\tau(\overline{\rho}_\tau(X))\}^{-1}\gamma_\tau\right).$$

For the gamma $G(dy; \tau) \propto y^{\alpha-1}\exp(-y/\beta)\,dy$, the MLE $\widehat{\tau} \equiv (\widehat{\alpha}, \widehat{\beta})$ satisfies

$$(3.7) \qquad \sum_{k=1}^{\infty} \frac{\sum_{\ell=k}^{\infty} n_\ell}{\widehat{\alpha} + k - 1} = \frac{\widetilde{d}\log(1+\widehat{\beta})}{1 - (1+\widehat{\beta})^{-\widehat{\alpha}}}, \qquad \frac{\widetilde{d}\widehat{\alpha}\widehat{\beta}}{1 - (1+\widehat{\beta})^{-\widehat{\alpha}}} = N,$$

with $\tilde{d} = \sum_{k=1}^{\infty} n_k$, and (3.4) holds [29]. Rao [19] called (3.3) with (3.7) pseudo MLE in a different (gamma) model, but the efficiency of the $\hat{d}$ was not clear [11].

The species problem is a special case of estimating (1.1) when $d$ is viewed as the number of species represented in the population out of a total of $n$ species. Specifically, letting $p_j = 0$ if the $j$th species is not represented in the population, estimating

$$(3.8) \qquad d = \sum_{j=1}^{n} I\{p_j > 0\} = \sum_{j=1}^{n} I\{X_j = 0, p_j > 0\} + \sum_{k=1}^{N} n_k$$

is equivalent to estimating (1.1) with $u(x,p) = I\{p > 0\}$ or $u(x,p) = I\{x = 0, p > 0\}$, based on observations $\{X_j, j \leq n\}$. Under (3.1) and (3.2) with $d$ replaced by $n$,

$$(3.9) \quad P_{p_*,\tau}\{X_j = k\} = (1-p_*)I\{k=0\} + p_* \frac{\int e^{-y}(y^k/k!)G_\tau(dy)}{\int (1-e^{-y})G_\tau(dy)} I\{k > 0\}$$

with certain $p_* < \int(1-e^{-y})G_\tau(dy)$. Under (3.9), the $\hat{\tau}$ in (3.3) is the conditional MLE of $\tau$ given $\{n_k, k \geq 1\}$. Since $(\sum_{k=1}^{\infty} n_k, d, n-d)$ is a trinomial vector, $\hat{\tau}$ in (3.3) equals the MLE of $\tau$ based on a sample $\{X_j, j \leq n\}$ from (3.9), provided that $\hat{d}$ in (3.3) is no greater than $n$. Since $P_{p_*,\tau}\{\hat{d} \leq n\} \to 1$ under (3.9), by Theorem 2.1, the (conditional) MLE (3.3) is asymptotically efficient in the empirical Bayes model (3.2) under conditions (2.2), (3.4) and (3.5).

3.2. *General mixture.* Now, suppose the distribution $G$ in (3.2) is completely unknown. The nonparametric MLE of $(d, G)$ is given by

$$(3.10) \quad \hat{d} \equiv \frac{\tilde{d} \int_{y>0} \hat{G}(dy)}{\int (1-e^{-y})\hat{G}(dy)}, \qquad \hat{G} \equiv \arg\max_G \prod_{k=1}^{\infty} \left\{ \frac{\int e^{-y} y^k G(dy)}{1 - \int e^{-y} G(dy)} \right\}^{n_k},$$

with $\tilde{d} \equiv \sum_{k=1}^{N} n_k$, but its asymptotic distribution is unclear. Since there is no solution $v$ to the equation $\sum_{x=0}^{\infty} v(x) e^{-\vartheta} \vartheta^x / x! = I\{\vartheta > 0\}$ for $0 \leq \vartheta < \infty$, by Theorems 2.3 and 2.5, the estimation of $d$ with completely unknown $G$ is an ill-posed problem.

Among many choices, a compromise between (3.3) and (3.10) is to fit $E_\tau n_k \propto P_\tau(X = k) = \int e^{-y}(y^k/k!) G_\tau(dy)$ for $1 \leq k \leq m$. For gamma $G$ with $En_{k+1}/En_k = (k+\alpha)\beta/(1+\beta)$, fitting the negative binomial distribution yields

$$(3.11) \qquad \hat{d} \equiv \tilde{d} + \max(\hat{\tau}_1, 0) n_1, \qquad \tilde{d} \equiv \sum_{k=1}^{N} n_k,$$

where $\widehat{\tau}_1$ is the (weighted) least squares estimate of $\tau_1 \equiv (\beta+1)/(\alpha\beta)$ based on

$$n_k = \tau_1 n_{k+1} + \tau_2 (k n_k) + \text{error}, \qquad k = 1, \ldots, m-1, \qquad \tau_2 \equiv -1/\alpha,$$

with $n_k$ being a response variable and $(n_{k+1}, k n_k)$ being covariates for each $k$. For small $\theta_j$ (large $n_k$ for small $k$), (3.11) has high efficiency for gamma $G$ and small bias for $G(y) = c_1 y^\alpha + (c_2 + o(1)) y^{\alpha+1}$ at $y \approx 0$. Chao [5] proposed $\widetilde{d} + n_1^2/(2n_2)$ as a low estimate of $d$. Another possibility is to estimate $d$ by correcting the bias of the estimator $\widetilde{d}/(1 - n_1/N)$ of Darroch and Ratcliff [9] as in [6].

**4. Networks: estimation of node degrees based on source-destination data.**
Source-destination (SD) data in networks are generated by sending probes (e.g., *traceroute* queries in the Internet) through networks from certain source nodes to certain destination nodes; see [8, 32]. We shall treat SD data as a collection of random vectors $W_j, j = 1, \ldots, N$, generated from a sample of SD pairs and make statistical inference based on $U$-processes of $\{W_j\}$, for example,

$$\text{(4.1)} \qquad \sum_{j=1}^{N} \frac{h_1(W_j)}{N}, \qquad \sum_{1 \leq j_1 \neq j_2 \leq N} \frac{h_2(W_{j_1}, W_{j_2})}{N(N-1)},$$

indexed by Borel $h_1$ and $h_2$, where $W_j$ are the observations from the $j$th SD pair in the sample. We focus here on the estimation of node degrees, although the approach based on (4.1) could be useful in other network problems.

The topology of a deterministic network can be described with a routing table: a list $r_1, \ldots, r_J$ of directed paths representing connections between pairs of source and destination nodes, with each path being composed of a set of directed links. For example, the path $4 \to 2 \to 3 \to 8$ has source node 4, destination node 8, and links $4 \to 2$, $2 \to 3$ and $3 \to 8$. Consider a network with nodes $\{1, \ldots, K\}$. The link degree $D(k, \ell)$ is defined as the number of paths using the link $k \to \ell$,

$$\text{(4.2)} \qquad D(k, \ell) \equiv \#\{j \leq J : \text{link } k \to \ell \text{ is used in } r_j\},$$

with $D(k, \ell) = 0$ if $k \to \ell$ is nonexistent or never used. The node degree, defined as

$$\text{(4.3)} \qquad d_k = \sum_{\ell=1}^{K} I\{D(k, \ell) > 0\},$$

is the number of outgoing links from $k$ to other nodes. This is also called out-degree. The in-degree, $\sum_\ell I\{D(\ell, k) > 0\}$, is the number of incoming links to $k$. The node degrees $d_k$ and their (empirical) distributions are important characteristics of networks; see [12, 15, 30].

For a given sample size $N$, let $R_1, \ldots, R_N$ be a sample of SD pairs from the routing table $\{r_1, \ldots, r_J\}$. Suppose we observe the paths of $R_j$, so that the vectors $W_j \equiv (W_{1j}, \ldots, W_{Kj})'$ are given by $W_{kj} \equiv \ell$ if link $k \to \ell$ is used in $R_j$ for some $1 \leq \ell \leq K$ and $W_{kj} = 0$ otherwise. The observed link frequencies are

$$(4.4) \quad X_{k\ell} \equiv \#\{j \leq N : \text{link } k \to \ell \text{ is used in } R_j\} = \sum_{j=1}^{N} I\{W_{kj} = \ell\}.$$

Since $X_{k\ell} = 0$ for $D(k, \ell) = 0$ by (4.3), the node degree $d_k$ is a sum

$$(4.5) \quad d_k = \widetilde{d}_k + s_k, \qquad \widetilde{d}_k \equiv \sum_{\ell=1}^{K} I\{X_{k\ell} > 0\},$$

where $\widetilde{d}_k$ is the observed degree and $s_k$ is the unobserved degree given by

$$(4.6) \quad s_k \equiv \sum_{\ell=1}^{K} I\{X_{k\ell} = 0, D(k, \ell) > 0\}.$$

Lakhina, Byers, Crovella and Xie [16] and Clauset and Moore [7] pointed out that the observed degrees $\widetilde{d}_k$ may grossly underestimate the true node degree $d_k$.

It follows from (4.5), (4.6) and (3.8) that the problem of estimating the node degree (4.3) is a species problem. From this point of view, we may directly use estimators in Section 3 and references therein, for example, (3.11). However, in network problems, we are typically interested in simultaneous estimation of many node degrees. Thus, information from $\{X_{k\ell}, \ell \leq K\}$ can be pooled from different nodes $k$. Let $\mathcal{K} \subseteq \{1, \ldots, K\}$ be a collection of "similar" and/or "independent" nodes. Let $\mathcal{G}$ be a family of distributions, for example, gamma with unit scale. Suppose the $G$ in (3.2) for different nodes are identical to a member of $\mathcal{G}$ up to scale parameters $\beta_k$. Then, as in (3.10), the (pseudo) MLE for $\{d_k, \beta_k, k \in \mathcal{K}, G\}$ is given by

$$(4.7) \quad \begin{aligned} \widehat{d}_k &\equiv \frac{\sum_{j=1}^{N} n_{kj} \int_{y>0} \widehat{G}(dy)}{\int (1 - e^{-\widehat{\beta}_k y}) \widehat{G}(dy)}, \\ (\widehat{\beta}, \widehat{G}) &\equiv \arg\max_{\beta, G} \prod_{k \in \mathcal{K}} \prod_{j=1}^{N} \left\{ \frac{\int e^{-\beta_k y} y^j G(dy)}{1 - \int e^{-\beta_k y} G(dy)} \right\}^{n_{kj}}, \end{aligned}$$

where $\beta \equiv (\beta, \ldots, \beta_K)$ and the maximum is taken over all $\beta_k > 0$ and $G \in \mathcal{G}$. This type of estimator is expected to perform well for self-similar networks.

In the nonparametric case of completely unknown $G$, the MLE $(\widehat{\beta}, \widehat{G})$ in (4.7) can be computed via the following EM algorithm:

$$\beta_k^{(m+1)} \leftarrow \left\{ \sum_{j=1}^{N} n_{kj} \left( \frac{p(j+1; \beta_k^{(m)}, G^{(m)})}{p(j; \beta_k^{(m)}, G^{(m)})} + \frac{p(1; \beta_k^{(m)}, G^{(m)})}{1 - p(0; \beta_k^{(m)}, G^{(m)})} \right) \right\}^{-1} \sum_{j=1}^{N} j n_{kj},$$

with $p(j; \beta_k, G) \equiv \int e^{-\beta_k y} y^j G(dy)$,

$$G^{(m+1)}(d\vartheta) \leftarrow G^{(m)}(d\vartheta) \left( \sum_{k \in \mathcal{K}} \sum_{j=1}^{N} n_{kj} / \{1 - p(0; \beta_k^{(m+1)}, G^{(m)})\} \right)^{-1}$$

$$\times \sum_{k \in \mathcal{K}} \sum_{j=1}^{N} n_{kj} \left( \frac{\exp(-\beta_k^{(m+1)} \vartheta) \vartheta^j}{p(j; \beta_k^{(m+1)}, G^{(m)})} + \frac{\exp(-\beta_k^{(m+1)} \vartheta)}{1 - p(0; \beta_k^{(m+1)}, G^{(m)})} \right).$$

**5. Data confidentiality: estimation of risk in statistical disclosure.** A major concern in releasing microdata sets is protecting the privacy of individuals in the sample. Consider a data set in the form of a high-dimensional contingency table. If an individual belongs to a cell with small frequency, an intruder with certain knowledge about the individual may identify him and learn sensitive information about him in the data. Statistical models and methods concerning the risk of such breach of confidentiality have been considered by many; see [10] and the proceedings of the joint ECE/EUROSTAT work sessions on statistical data confidentiality. For multi-way contingency tables, Polettini and Seri [18] and Rinott [21] studied the estimation of global disclosure risks of the form

$$(5.1) \qquad S_J \equiv \sum_{j=1}^{J} u(X_j, Y_j)$$

based on $\{X_j, j \leq J\}$, where $X_j$ and $Y_j$ are the sample and population frequencies in the $j$th cell, $J$ is the total number of cells, and $u(x, y)$ is a loss function of the form $u(x, y) = u(x)/y$, for example, $u(x, y) = y^{-1} I\{x = 1\}$.

Let $N = \sum_{j=1}^{J} Y_j$ be the population size. Suppose $N \sim \text{Poisson}(\lambda)$,

$$(5.2) \quad \{Y_j\} | N \sim \text{multinomial}(N, \{\pi_j\}), \qquad X_j | (\{Y_j\}, N) \sim \text{binomial}(Y_j, p_j),$$

for certain $\pi_j > 0$ with $\sum_{j=1}^{J} \pi_j = 1$, $0 \leq p_j \leq 1$ and $\lambda > 0$. For known $\{p_j, \pi_j, \lambda\}$, the Bayes estimator of $S_J$ in (5.1) is

$$(5.3) \quad S_J^* \equiv E(S_J | \{X_j\}) = \sum_{j=1}^{J} \overline{u}_j(X_j), \qquad \overline{u}_j(x) \equiv E u(x, Y_j - X_j + x),$$

with $Y_j - X_j \sim \text{Poisson}((1 - p_j) \pi_j \lambda)$ (independent of $X_j$). For $u(x, y) = y^{-1} I\{x = 1\}$,

$$(5.4) \qquad \overline{u}_j(x) = \{(1 - p_j) \pi_j \lambda\}^{-1} [1 - \exp\{-(1 - p_j) \pi_j \lambda\}].$$

In general, the parameters $(1 - p_j) \pi_j \lambda$ cannot be completely identified from the data $X_j \sim \text{Poisson}(p_j \pi_j \lambda)$, so that it is necessary to further model the parameters. This can be achieved by setting $\{p_j, \pi_j, \lambda\}$ to known tractable functions of an unknown vector $\tau$ and certain covariates $z_j$ characterizing cells $j$, and by incorporating all available knowledge about the parameters, for example, $\lambda \approx N$ and $\sum_{j=1}^{J} p_i \pi_j \approx n/N$, where $n = \sum_{j=1}^{J} X_j$ is the sample size. Consequently, the conditional expectation $\overline{u}_j(x)$ in (5.4) can be written as $\overline{u}_j(x) = \overline{u}(x, z_j; \tau)$. This suggests

$$(5.5) \qquad \widehat{S}_J \equiv \sum_{j=1}^{J} \overline{u}(X_j, z_j; \widehat{\tau}_J)$$

as an estimator of the global risk (5.1) and its conditional expectation (5.3), where $\widehat{\tau}_J$ is a suitable (e.g., the maximum likelihood or method of moments) estimator of $\tau$. For example, in a two-way table with cells labelled by $j \sim (i,k)$ and known $\pi_{i,k}$ and $\lambda$, we may assume a regression model $p_{i,k} = \psi_0(\tau_1 + \tau_2' z_{i,k})$ for a certain known (e.g., logit or probit) function $\psi_0$. In the case of unknown $\pi_{i,k}$, we may consider the independence model $\pi_{i,k} = \pi_i.\pi_{\cdot k}$ with unknown $\pi_i.$ and known or unknown $\pi_{\cdot k}$. If $\tau$ has fixed dimensionality and $\widehat{\tau}_J$ is asymptotically efficient, (5.5) is efficient by Theorem 2.2. Theorem 2.2 also suggests that (5.5) is highly efficient if $\dim(\tau)/J \to 0$.

Alternatively, we may consider the negative binomial model $N \sim \text{NB}(\alpha, 1/(1+\beta))$, that is, $P(N=k) = \Gamma(k+\alpha)\{\Gamma(\alpha)k!\}^{-1}\beta^k/(1+\beta)^{k+\alpha}$. As in [21], we have in this case $Y_j \sim \text{NB}(\alpha, 1/(1+\beta_j))$ with $\beta_j = \beta\pi_j$, $X_j \sim \text{NB}(\alpha, 1/(1+p_j\beta_j))$, and $(Y_j - X_j)|\{X_j = x\} \sim \text{NB}(x+\alpha, (1+p_j\beta_j)/(1+\beta_j))$. Consequently,

$$(5.6) \qquad \overline{u}_j(x) = \frac{1+p_j\beta_j}{(1-p_j)\beta_j} \int_{(1+p_j\beta_j)/(1+\beta_j)}^{1} t^{\alpha_j - 1} dt\, I\{x=1\}$$

in (5.3) for $u(x,y) = y^{-1} I\{x=1\}$. Bethlehem, Keller and Pannekoek [2] studied this negative binomial model with constant $\pi_j = 1/J$ and $p_j = En/EN \approx n/N$. For $(\alpha_j, \beta_j) \to (0, \infty)$, $(Y_j - X_j)|\{X_j = x\}$ converges in distribution to the $\text{NB}(x, p_j)$, resulting in the $\mu$-ARGUS estimator [1] with $\overline{u}_j(x) = p_j(1-p_j)^{-1}(-\log p_j)I\{x=1\}$ in (5.6), as pointed out by Rinott [21]. Compared with the Poisson model in which $\lambda \approx N$, estimates of both $EN$ and $\text{Var}(N)$ are required in the negative binomial model. The $\mu$-ARGUS model essentially assumes $\text{Var}(N)/(EN)^2 \geq 1/\alpha \to \infty$, which may not be suitable in some applications.

**6. General information bounds.** We provide a lower bound for the asymptotic variance and a convolution theorem for (locally asymptotically) regular estimators of the sum in (1.2). To facilitate the statements of our results, we first briefly describe certain terminologies and concepts in general asymptotic theory.

6.1. *Scores and tangent spaces.* Suppose $(X, \theta) \sim F$ with $F \in \mathcal{F}$, where $\mathcal{F}$ is a family of joint distributions. Let $\mathcal{C} \equiv \mathcal{C}(F_0)$ be a collection of mappings $\{F_t, 0 \leq t \leq 1\}$ from $[0,1]$ to $\mathcal{F}$ satisfying

(6.1)  $E_{F_0}(\sqrt{f_t(X)} - 1 - t\rho(X)/2)^2 = o(t^2), \qquad E_{F_0} f_t(X) = 1 + o(t^2),$

for certain score functions $\rho(x) \equiv \rho(x; \{F_t\})$ depending on the mappings $\{F_t\}$, where $f_t \equiv dF_t^X/dF_0^X$ is the Radon–Nikodym derivative of the absolutely continuous part of the marginal distribution $F_t^X$ of $X$ under $F_t$ with respect to the marginal distribution $F_0^X$. Let $\mathcal{C}_* \equiv \mathcal{C}_*(F_0)$ be the collection of score functions $\rho(X)$ generated by $\mathcal{C}$. The tangent space $H_* \equiv H_*(F_0)$ is the closure of the linear span $[\mathcal{C}_*]$ of $\mathcal{C}_*$ in $L_2(F_0)$; that is,

(6.2)  $\qquad H_* \equiv \overline{[\mathcal{C}_*]}, \qquad \mathcal{C}_* \equiv \{\rho(\cdot; \{F_t\}) : \{F_t\} \in \mathcal{C}\}.$

For further discussion about score and tangent space, see [3], pages 48–57. The second part of (6.1) holds in regular parametric models; see [3], page 459.

6.2. *Smoothness of random variables and their distributions.* Let $\mathcal{L}(U; F)$ be the distribution of $U$ under $P_F$. Suppose that, for all $\{F_t\} \in \mathcal{C}$, the random variables $u_{F_t} \equiv u(X, \theta; F_t)$ and $\overline{u}_{F_t} \equiv E_{F_t}[u_{F_t}|X]$ satisfy the continuity conditions

(6.3)  $\qquad \lim_{t \to 0+} \mathrm{Var}_{F_0}(\overline{u}_{F_t} - \overline{u}_{F_0}) = 0,$

(6.4)  $\qquad \mathcal{L}(w_{F_t}; F_t) \xrightarrow{D} \mathcal{L}(w_{F_0}; F_0), \qquad E_{F_t} w_{F_t}^2 \to E_{F_0} w_{F_0}^2,$

as $t \to 0+$, with $w_F \equiv \overline{u}_F - u_F$, and also satisfy the differentiability condition

(6.5)  $\qquad \lim_{t \to 0+} E_{F_0}(\overline{u}_{F_t} - \overline{u}_{F_0})/t = E_{F_0} \phi(X) \rho(X)$

for certain $\phi(X) \equiv \phi(X; F_0) \in L_2(F_0)$. The usual smoothness condition for $\mu(F)$, see [3], pages 57–58, is that, for a certain influence function $\psi(X) \equiv \psi(X; F_0) \in L_2(F_0)$,

(6.6)  $\qquad \lim_{t \to 0+} \{\mu(F_t) - \mu(F_0)\}/t = E_{F_0} \psi(X) \rho(X).$

6.3. *Regular estimators.* An estimator $\tilde{\mu}_n \equiv \tilde{\mu}_n(X_1, \ldots, X_n)$ of $\mu(F)$ is (locally asymptotically) regular at $F_0$ if there exists a random variable $\zeta_0$ such that

(6.7)  $\qquad \lim_{n \to \infty} \mathcal{L}(n^{1/2}\{\tilde{\mu}_n - \mu(F_{c/\sqrt{n}})\}; F_{c/\sqrt{n}}) = \mathcal{L}(\zeta_0; F_0)$

for all $c > 0$ and $\{F_t\} \in \mathcal{C}$ ([3], page 21). Likewise, for the estimation of the sum $S_n(F)$ in (1.2), we say that an estimator $\widetilde{S}_n \equiv \widetilde{S}_n(X_1, \ldots, X_n)$ is regular at $F_0$ if there exists a random variable $\xi_0$ such that, for all $c > 0$ and $\{F_t\} \in \mathcal{C}$,

(6.8)  $\qquad \lim_{n \to \infty} \mathcal{L}(n^{-1/2}\{\widetilde{S}_n - S_n(F_{c/\sqrt{n}})\}; F_{c/\sqrt{n}}) = \mathcal{L}(\xi_0; F_0).$

6.4. *Efficient influence functions and information bounds.* Let $\psi_*$ be the projection of $\psi$ in (6.6) to the tangent space $H_*$ in (6.2). The standard convolution theorem ([3], page 63) asserts that, for a certain variable $\zeta_0'$,

$$\mathcal{L}(\zeta_0; F_0) = N(0, E\psi_*^2(X)) \star \mathcal{L}(\zeta_0'; F_0)$$

for the $\zeta_0$ in (6.7), and that efficient estimators are characterized by (1.4). For $h \in L_2(F_0)$, let $A_n(h) \equiv \sum_{j=1}^n h(X_j, \theta_j)/n$ and $Z_n(h) \equiv \sqrt{n}\{A_n(h) - E_{F_0}h\}$.

THEOREM 6.1. *Suppose* (6.3), (6.4) *and* (6.5) *hold at* $F_0$. *Let* $\phi_{*,0}$ *be the projection of* $\phi$ *in* (6.5) *into the tangent space* $H_*$ *in* (6.2), *and let* $\phi_* \equiv \overline{u}_{F_0} - \mu(F_0) + \phi_{*,0}$.

(i) *If* (6.8) *holds, then* $\mathrm{Var}_{F_0}(\xi_0) \geq \mathrm{Var}_{F_0}(\phi_* - u_{F_0})$. *Moreover, the lower bound is reached without bias, that is,* $E_{F_0}\xi_0^2 = \mathrm{Var}_{F_0}(\phi_* - u_{F_0})$, *iff* (1.5) *holds.*

(ii) *If* (6.8) *holds and the* $L_2(F_0)$ *closure* $\overline{\mathcal{C}}_*$ *of* $\mathcal{C}_*$ *in* (6.2) *is convex, then there exist a random variable* $\tilde{\xi}_0$ *and certain normal variables* $Z(h) \sim N(0, \mathrm{Var}_{F_0}(h))$ *such that*

$$\mathcal{L}\left(\begin{pmatrix} \sqrt{n}\{\widetilde{S}_n/n - A_n(\phi_*) - \mu(F_0)\} \\ Z_n(\overline{u}_{F_0} + h - u_{F_0}) \end{pmatrix}; F_0\right) \xrightarrow{D} \mathcal{L}\left(\begin{pmatrix} \tilde{\xi}_0 \\ Z(\overline{u}_{F_0} + h - u_{F_0}) \end{pmatrix}; F_0\right)$$

*and* $\tilde{\xi}_0$ *is independent of* $Z(\overline{u}_{F_0} + h - u_{F_0})$ *for all* $h \in H_*$. *In particular, for* $h = \phi_{*,0}$,

$$\mathcal{L}(\xi_0; F_0) = \mathcal{L}(Z(\phi_* - u_{F_0}); F_0) \star \mathcal{L}(\tilde{\xi}_0; F_0).$$

(iii) *Suppose* $E_{F_t}\overline{u}^2(X; F_t)$ *is bounded for all* $\{F_t\} \in \mathcal{C}$. *Then,* $\psi_* = \phi_{*,0} + \overline{u}_*$ *is the efficient influence function for the estimation of* $\mu(F)$, *that is,* (6.6) *holds with* $\psi = \psi_*$, *where* $\overline{u}_*$ *is the projection of* $\overline{u}_{F_0}$ *to* $H_*$. *Consequently,* (1.6) *holds.*

REMARK 6.1. Based on Theorem 6.1(i) and (ii), $\widehat{S}_n$ is said to be locally asymptotically efficient if (1.5) holds. Note that in Theorem 6.1(ii), $\tilde{\xi}_0 = 0$ iff (1.5) holds.

REMARK 6.2. In the proof of Theorem 6.1(iii), we show that (6.5) and (6.6) are equivalent under the condition that $E_{F_t}\overline{u}^2(X; F_t) = O(1)$ for all $\{F_t\} \in \mathcal{C}$.

REMARK 6.3. For the estimation of $\mu(F)$, that is, $u(x, \vartheta, F) \equiv \mu(F)$ as a special case of Theorem 6.1(ii), a standard proof of the convolution theorem uses analytic continuation along lines passing through the origin in the tangent space, and as a result, $\overline{\mathcal{C}}_*$ is often assumed to be a linear space. In the proof of Theorem 6.1(ii), analytic continuation is used along arbitrary lines across $\overline{\mathcal{C}}_*$, so that only the convexity of $\overline{\mathcal{C}}_*$ is needed as in [31], pages 366–367. Rieder [20] showed that, in the case of convex $\overline{\mathcal{C}}_*$, the projections of scores to $\overline{\mathcal{C}}_*$ (not to $H_*$) are useful in the context of one-sided confidence.

6.5. *Finite-dimensional models.* Let $\mathcal{F} = \{F_\tau, \tau \in \mathcal{T}\}$ with an open Euclidean parameter space $\mathcal{T}$. We shall extend the results in Section 2.1 to general sums (1.2). Suppose $dF_\tau^X = f_\tau^X d\nu$ exists and is differentiable in the sense of (6.1), that is,

$$(6.9) \qquad \int (f_{\tau+\Delta}^{1/2} - f_\tau^{1/2} - \Delta \rho_\tau)^2 d\nu = o(\|\Delta\|^2), \qquad \tau \in \mathcal{T}.$$

Let $E_\tau \equiv E_{F_\tau}$, $I_\tau \equiv \mathrm{Cov}_\tau(\rho_\tau(X))$, $u_\tau \equiv u(X, \theta; F_\tau)$ and $\overline{u}_\tau \equiv \overline{u}(X; F_\tau)$.

THEOREM 6.2. (i) *Suppose* (6.9) *holds*, $I_\tau$ *is of full-rank*, $\mathcal{L}(u_\tau; F_\tau)$ *is continuous in $\tau$ in the weak topology,* $E_\tau u_\tau^2$ *is continuous,* $E_\tau \{\overline{u}_{\tau+\Delta} - \overline{u}_\tau\}^2 \to 0$ *as* $\Delta \to 0$, $E_\tau \overline{u}_\tau^2$ *is locally bounded, and* $\mu'(\tau)$ *exists. Then* (2.4) *gives the efficient influence function for the estimation of* (1.2) *with* $\gamma_\tau = \mu'(\tau) - E_\tau \overline{u}_\tau \rho_\tau$, *and* (1.5) *and* (1.6) *hold.*

(ii) *Suppose* (2.6), (2.7) *and conditions of* (i) *hold. Then* (2.8) *holds for the plug-in estimator* (2.5) *with the $\gamma_\tau$ in* (i). *In particular,* (2.5) *is asymptotically efficient under* $P_\tau$ *iff* $\gamma_\tau \kappa_\tau = \gamma_\tau I_\tau^{-1} \rho_\tau$.

REMARK 6.4. Comparing Theorem 6.2 with Theorems 2.1 and 2.2, we see that (6.9) is weaker than (2.2) and (1.2) is more general than (1.1), while stronger conditions are imposed on $u_\tau$ in Theorem 6.2.

**7. Proofs.** We prove Theorems 6.1, 2.1, 2.2, 6.2, and 2.3–2.5 in this section.

LEMMA 7.1. *Suppose* (2.2) *holds. Let* $(X, \theta) \sim F_t$ *under* $P_{\tau+at}$ *and* $\rho = a^t \rho_\tau$ *for a vector $a$, where $\rho_\tau$ is as in* (2.3). *Then* (6.1) *holds with* $P_{F_0} = P_\tau$.

PROOF. Let $g_t \equiv g_{\tau+at}$ and $\Delta = at$. The lemma follows from the expansion

$$\frac{\sqrt{f_t} - 1}{t} - \frac{\rho}{2} = \frac{1}{f_t^{1/2} + 1} E_0 \left[ \frac{g_t^{1/2} - 1}{t} (g_t^{1/2} + 1) \Big| X = x \right] - E_0 \left[ \frac{a^t \widetilde{\rho}_\tau}{2} \Big| X = x \right].$$

The uniform integrability of the square of the right-hand side (i.e., the first term) under $f_0(x)$ follows from the inequality $E_0[g_t|X] \le f_t(X)I\{f_0(X) > 0\}$. We omit the details. □

LEMMA 7.2. *Suppose* (6.1) *holds and* $X \sim F_t^X$ *under* $P_t$, $0 \leq t \leq 1$. *Let* $\mu_t \equiv E_t h_t(X)$ *for a certain Borel* $h_t$. *If* $E_t h_t^2(X) = O(1)$ *and* $h_t \to h_0$ *in* $L_2(P_0)$, *then*

$$\mu_t - \mu_0 = E_0\{h_t(X) - h_0(X)\} + tE_0\rho(X)h_0(X) + o(t) \qquad \text{as } t \to 0.$$

PROOF. Let $B_t$ be the support sets of $dP_t(X) - f_t(X)\,dP_0(X)$. By (6.1) and the boundedness of $E_t h_t^2$, $E_t h_t - E_0 f_t h_t = E_t h_t I_{B_t} = O(1)(E_t h_t^2)^{1/2} \times P_t^{1/2}(B_t) = o(t)$. Thus,

(7.1) $\quad \mu_t - \mu_0 = E_t h_t - E_0 h_0 = E_0(f_t - 1)h_t + E_0(h_t - h_0) + o(t)$

as $t \to 0+$. Since $(\sqrt{f_t} - 1)/t \to \rho/2$ in $L_2(P_0)$ and $E_0\{(\sqrt{f_t}+1)h_t\}^2 = O(1)$,

$$E_0(f_t - 1)h_t/t = E_0[t^{-1}(\sqrt{f_t} - 1)(\sqrt{f_t} + 1)h_t] \to E_0 h_0 \rho.$$

This and (7.1) complete the proof. $\square$

PROOF OF THEOREM 6.1. Let $F_n \equiv F_{c/\sqrt{n}}$, $\xi_n \equiv \sqrt{n}\{\widetilde{S}_n/n - S_n(F_n)/n\}$, $\xi_n' \equiv \sqrt{n}\{\widetilde{S}_n/n - A_n(\overline{u}_{F_n})\}$, $\xi_n'' \equiv \sqrt{n}A_n(w_{F_n})$ and $Z'' = Z(w_{F_0})$. Then $\xi_n = \xi_n' + \xi_n''$ and $\xi_n'$ depend on $\{X_j\}$ only. By (6.4), $w_{F_n}^2$ under $P_{F_n}$ are uniformly integrable and $\mathcal{L}(w_{F_n}; F_n) \xrightarrow{D} \mathcal{L}(w_{F_0}; F_0)$ as $n \to \infty$. Thus, by the Lindeberg central limit theorem and the weak law of large numbers,

(7.2) $\qquad E_{F_n}[\exp(it\xi_n'')|\{X_j\}] \to E_{F_0}\exp(itZ'')$

in probability for all $t$. Since $\xi_n'$ depends on $\{X_j\}$ only, this and (6.8) imply

$$E_{F_n}\exp(it\xi_n')E\exp(itZ'') = E_{F_n}\exp(it\xi_n')\exp(it\xi_n'') + o(1) \to E_{F_0}\exp(it\xi_0).$$

Thus, since $E\exp(itZ'') \neq 0$ for all $t$,

(7.3) $\quad \mathcal{L}\left(n^{-1/2}\left\{\widetilde{S}_n - \sum_{j=1}^n \overline{u}(X_j; F_{c/\sqrt{n}})\right\}; F_{c/\sqrt{n}}\right) = \mathcal{L}(\xi_n'; F_n) \xrightarrow{D} \mathcal{L}(\xi_0'; F_0)$

for a certain variable $\xi_0'$ independent of $c > 0$ and the curve $\{F_t\} \in \mathcal{C}$.

Define $\xi_{n,0}' \equiv \sqrt{n}\{\widetilde{S}_n/n - A_n(\overline{u}_{F_0})\}$. By (6.3) and (6.5), $\xi_{n,0}' - \xi_n' = \sqrt{n}A_n \times (\overline{u}_{F_n} - \overline{u}_{F_0}) = E_{F_0}(\overline{u}_{F_n} - \overline{u}_{F_0}) + o_P(1) \to cE\phi(X)\rho(X)$ in probability under $P_{F_0}$. Thus, as in [3], pages 24–26, by (7.3) and the LAN from (6.1) and (6.2),

(7.4) $\quad E_{F_0}\exp(it\xi_0' + zZ(\rho)) = \exp[itzE_{F_0}\phi\rho + z^2 E_{F_0}\rho^2/2]E_{F_0}\exp(it\xi_0')$

for all $\rho \in \mathcal{C}_*$ and complex $z$. Here $Z(h)$ are constructed so that $(\xi_{n,0}', Z_n(h))$ converges jointly in distribution to $(\xi_0', Z(h))$ for all $h \in L_2(F_0)$. Differentiating (7.4) in $t$ at $t = 0$ and then in $z$ at $z = 0$, we find

(7.5) $\qquad E_{F_0}\xi_0' Z(h) = E_{F_0}\phi(X)h(X) = E_{F_0}Z(\phi_{*,0})Z(h)$

for all scores $h = \rho$, $\rho \in \mathcal{C}_*$, and then for all $h \in H_*$ by (6.2). Since $\phi_{*,0} \in H_*$, $\xi'_0 - Z(\phi_{*,0})$ and $Z(\phi_{*,0})$ are orthogonal in $L_2(F_0)$. This proves (i), since $\xi'_0$ and $Z(\phi_{*,0})$ are both independent of $Z''$ by (7.2) and $Z(\phi_{*,0}) + Z'' = Z(\phi_* - u_{F_0})$.

Now, suppose $\overline{\mathcal{C}_*}$ is convex in $L_2(F_0)$. By continuity extension, (7.4) holds for all $\rho \in \overline{\mathcal{C}_*}$ and complex $z$. Let $\rho_j \in \overline{\mathcal{C}_*}$. Since (7.4) holds for $\rho = s\rho_1 + (1-s)\rho_2, 0 \leq s \leq 1$, with both sides being analytic in $s$, by analytic continuation it holds for $\rho = s\rho_1 + (1-s)\rho_2$ for all real $s$. Thus, (7.4) holds for all complex $z$ and

$$(7.6) \qquad \rho \in H_0 \equiv \{s\rho_1 + (1-s)\rho_2 : \rho_j \in \overline{\mathcal{C}_*}, -\infty < s < \infty\}.$$

Let $\widetilde{H}$ be the linear span of a set of finitely many members of $\overline{\mathcal{C}_*}$. Let $\rho_1$ be a fixed interior point of $\widetilde{H} \cap \overline{\mathcal{C}_*}$ and $\rho_2 \in \widetilde{H}$ with $\|\rho_2 - \rho_1\| = \delta_0$. For sufficiently small $\delta_0 > 0$, $\rho_2 \in \overline{\mathcal{C}_*}$ for all such $\rho_2$, so that $\widetilde{H} \subseteq H_0$. Thus, $H_0$ is a linear space and $H_*$ is the closure of $H_0$. It follows that (7.4) holds for all $\rho \in H_*$ and complex $z$. As in [3], pages 25–26, this implies the independence of $\xi'_0 - Z(\phi_{*,0})$ and $\{Z(h) : h \in H_*\}$. Since $\{\xi'_0, Z(h), h \in H_*\}$ is independent of $Z'' = Z(\overline{u}_{F_0} - u_{F_0})$ by (7.2), the conclusions of part (ii) hold with $\tilde{\xi}_0 = \xi'_0 - Z(\psi_{*,0})$.

The proof of part (iii) follows easily from Lemma 7.2 with $h_t = \overline{u}_{F_t}$, which gives

$$\{\mu(F_t) - \mu(F_0)\}/t - E_{F_0}\{\overline{u}_{F_t} - \overline{u}_{F_0}\}/t \to E_{F_0}\overline{u}_{F_0}\rho = E_{F_0}\overline{u}_*\rho.$$

It follows that (6.5) and (6.6) are equivalent under $E_{F_t}\overline{u}^2(X; F_t) = O(1)$, with $\psi = \psi_* = u_* + \phi_{*,0}$, by (1.6) and the definition of $\phi_*$. The proof is complete. □

PROOF OF THEOREM 2.1. The proof is similar to that of Theorem 6.1(i), so we omit certain details. By (2.2), $\xi_0$ is independent of $Z(\widetilde{\rho}_\tau)$ under $P_\tau$. Since $E_\tau u^2 < \infty$, (7.2) holds for fixed $F_n = F_\tau$, so that $\xi_0 = \xi'_0 + Z(\overline{u}_\tau - u)$ as a sum of independent variables. Let $Z(h_\tau)$ be the projection of $\xi'_0$ to $\{Z(h), h \in L_2(F_\tau)\}$ in $L_2(P_\tau)$ and $v_\tau = h_\tau + \overline{u}_\tau$. Then $\mathrm{Var}_\tau(\xi_0) \geq E_\tau(v_\tau - u)^2$ and $E_\tau(v_\tau - u)\widetilde{\rho}_\tau = 0$. Since $\xi'_0$ is the limit of variables dependent on $\{X_j\}$ only, $h_\tau$ and $v_\tau$ depend on $X$ only.

Since $E_\tau u^2 g_{\tau,\Delta}(\theta) \leq E_{\tau+\Delta} u^2 = O(1)$, by (2.2) and Lemma 7.2 with $h_t = h_0 = u(x, \vartheta)$, $\mu_{\tau+\Delta} - \mu_\tau \approx \Delta^t E_\tau u\widetilde{\rho}_\tau = \Delta^t E_\tau \psi_{*,\tau}(X)\rho_\tau(X)$, where $\psi_{*,\tau} \equiv \rho_\tau^t I_\tau^{-1} E_\tau u\widetilde{\rho}_\tau$. It follows that $0 = E_\tau(v_\tau - u)\widetilde{\rho}_\tau = E_\tau(v_\tau\widetilde{\rho}_\tau - \psi_{*,\tau}\rho_\tau) = E_\tau(v_\tau - \psi_{*,\tau})\rho_\tau$. Thus, $E_\tau(v_\tau - \overline{u}_\tau)\rho_\tau = E_\tau(\psi_{*,\tau} - \overline{u}_{*,\tau})\rho_\tau$ with $\overline{u}_{*,\tau} \equiv \rho_\tau^t I_\tau^{-1} E_\tau \overline{u}_\tau \rho_\tau$. Since $\psi_{*,\tau} - \overline{u}_{*,\tau}$ is linear in $\rho_\tau$, $Z(v_\tau - \overline{u}_\tau - (\psi_{*,\tau} - \overline{u}_{*,\tau}))$ is independent of $Z(\psi_{*,\tau} - \overline{u}_{*,\tau})$. Thus, $\mathrm{Var}_\tau(v_\tau - \overline{u}_\tau) \geq \mathrm{Var}_\tau(\psi_{*,\tau} - \overline{u}_{*,\tau}))$ and $\mathrm{Var}_\tau(\xi_0) \geq \mathrm{Var}_\tau(v_\tau - \overline{u}_\tau) + \mathrm{Var}_\tau(\overline{u}_\tau - u) \geq \mathrm{Var}_\tau(\phi_{*,\tau} - u)$ by (2.4). The proof is complete. □

PROOFS OF THEOREMS 2.2 AND 6.2. Theorem 6.2(i) follows from Theorem 6.1 and Remark 6.2. Let $\mu(t;\tau) = E_\tau \overline{u}_t(X)$. By Lemma 7.2, $\mu' = E_\tau u \widetilde{\rho}$ in Theorem 2.2 and $\gamma_\tau = (\partial/\partial t)\mu(\tau;\tau)$ in both theorems. Simple expansion of (2.5) via (2.7) yields

$$\frac{\widehat{S}_n}{n} = A_n(\overline{u}_\tau) + \{\mu(\widehat{\tau}_n;\tau) - \mu(\tau;\tau)\} + o_{P_\tau}(n^{-1/2})$$
$$= A_n(\overline{u}_\tau + \gamma_\tau \kappa_\tau) + o_{P_\tau}(n^{-1/2}),$$

which implies (2.8). Note that $\gamma_\tau(\kappa_\tau - \kappa_{*,\tau})$ is orthogonal to $\overline{u}_\tau - u_\tau + \gamma_\tau \kappa_{*,\tau}$. The proof is complete. $\square$

PROOFS OF THEOREMS 2.3, 2.4 AND 2.5. Let $G_t \equiv (1-t)G_0 + tG$, $f_t \equiv f_{G_t}$ and $E_t \equiv E_{G_t}$, $t > 0$. By (2.9), (6.1) holds with $\rho = f_G/f_0 - 1$. Since $E_G u^2 < \infty$, $u^2$ are uniformly integrable under $P_t$, so that (6.4) holds. Since $f_0/f_t \leq 1/(1-t)$, $\{\overline{u}_t^2, 0 \leq t \leq 1/2\}$ are uniformly integrable under $E_0$, so that (6.3) holds. Moreover,

$$(7.7) \quad t^{-1} E_0\{\overline{u}_t - \overline{u}_0\} = E_0\left\{\frac{f_G}{f_t}(\overline{u}_G - \overline{u}_0)\right\} \to E_0\left\{\frac{f_G}{f_0}(\overline{u}_G - \overline{u}_0)\right\}.$$

Suppose there exists a regular estimator of (1.1). Let $\xi_0'$ be as in (7.5) and let $Z(v - \overline{u}_0)$ be the projection of $\xi_0'$ to $\{Z(h), h \in L_2(f_0)\}$ as in the proof of Theorem 2.1. It follows from (7.7) and the argument leading to (7.5) that

$$E_0(v - \overline{u}_0)(f_G/f_0 - 1) = E_0 Z(v - \overline{u}_0) Z(\rho) = E_0\left\{\frac{f_G}{f_0}(\overline{u}_G - \overline{u}_0)\right\},$$

which implies $E_G v - E_0 v + E_0 u = E_G u$. Since $\xi_0'$ does not depend on the choice of $G \in \mathcal{G}_{G_0}$, $v \in \mathcal{V}_{G_0}$. By the Lindeberg central limit theorem, $E_{G_0} v^2 < \infty$ and $v \in \mathcal{V}_{G_0}$ imply $\mathcal{L}(Z_n(v-u); P_{c/\sqrt{n}}) \to \mathcal{L}(Z(v-u); P_0)$, so that $V_n$ in (2.15) is regular at $G_0$ for all $v \in \mathcal{V}_{G_0}$. If $v$ is a limit point of $\mathcal{V}_{G_0}$ in $L_2(f_0)$, $V_n$ is also a regular estimator of $S_n$ at $P_0$, so that $\mathcal{V}_{G_0}$ is closed in $L_2(f_0)$. This completes the proof of Theorem 2.3.

The proof of Theorem 2.4 is similar to those of Theorems 2.2 and 6.2 but simpler. We note that $E_{G_0}(v_G - v_{G_0}) = 0$. Finally, Theorem 2.5 follows from the fact that $\mathcal{V}_G$ contains a single function $v$ due to the completeness of exponential families. The proofs are complete. $\square$

DEPARTMENT OF STATISTICS
RUTGERS UNIVERSITY
HILL CENTER
BUSCH CAMPUS
PISCATAWAY, NEW JERSEY 08854-8019
USA
E-MAIL: czhang@stat.rutgers.edu